\documentclass[leqno,11pt]{amsart}
\usepackage{graphicx}
\usepackage{amscd}
\usepackage{amsmath}
\usepackage{caption}
\usepackage{amsfonts}
\usepackage{amssymb}
\usepackage{mathrsfs}
\usepackage{multicol}
\usepackage{color}
\usepackage{mathdots}
\usepackage{bigints}
\usepackage[labelsep=newline,font={small,up},margin=0pc]{caption}
\newcommand{\R}{\mathbb{R}}

\usepackage{tikz}
\usepackage{tikz-3dplot}
\usetikzlibrary{arrows,decorations,patterns,shadings,calc,shapes}
\usepackage{pgfplots}
\pgfplotsset{compat=1.12}

\newcommand{\be}{\begin{equation}}
\newcommand{\ee}{\end{equation}}
\newcommand{\ben}{\begin{eqnarray*}}
\newcommand{\een}{\end{eqnarray*}}

\newtheorem{thm}{Theorem}
\newtheorem{lem}{Lemma}
\newtheorem{cor}{Corollary}
\newtheorem{Prop}{Proposition}

\theoremstyle{definition}
\newtheorem{rem}{Remark}

\allowdisplaybreaks
\definecolor{darkgreen}{rgb}{0.09, 0.45, 0.27}
\definecolor{debianred}{rgb}{0.84, 0.04, 0.33}
\allowdisplaybreaks[1]
\newcommand{\divergence}{\mathop{\mathrm{div}}\nolimits}

\newcommand{\eps}{\varepsilon}
\newcommand{\dist}{\mathop{\mathrm{dist}}\nolimits}
\newcommand{\supp}{\mathop{\mathrm{supp}}\nolimits}
\usepackage{stackrel}
\usepackage{url}
\usepackage[normalem]{ulem}
\definecolor{DarkGreen}{rgb}{0,0.5,0.1} 

\newcommand\soutD{\bgroup\markoverwith
{\textcolor{DarkGreen}{\rule[.5ex]{2pt}{1pt}}}\ULon}
\newcommand\soutP{\bgroup\markoverwith
{\textcolor{blue}{\rule[.5ex]{2pt}{1pt}}}\ULon}
\newcommand{\Hm}[1]{\leavevmode{\marginpar{\tiny%
$\hbox to 0mm{\hspace*{-0.5mm}$\leftarrow$\hss}%
\vcenter{\vrule depth 0.1mm height 0.1mm width \the\marginparwidth}%
\hbox to
0mm{\hss$\rightarrow$\hspace*{-0.5mm}}$\\\relax\raggedright #1}}}

\begin{document}

\title[Curved nonlinear waveguides]
{Curved nonlinear waveguides} 

\author[L.~Baldelli]{Laura Baldelli} \email{labaldelli@ugr.es}

\author[D.~ Krej\v ci\v r\'ik]{David  Krej\v ci\v r\'ik}
\email{david.krejcirik@fjfi.cvut.cz}

\address[Baldelli]{IMAG, Departamento de An\'alisis Matem\'atico, 
Universidad de Granada, Campus Fuentenueva, 18071 Granada, Spain}
\address[ Krej\v ci\v r\'ik]{Department of Mathematics, 
Faculty of Nuclear Sciences and Physical Engineering, 
Czech Technical University in Prague, 
Trojanova 13, 12000 Prague, Czech Republic}


\begin{abstract}
The Dirichlet $p$-Laplacian in tubes of arbitrary cross-section 
along infinite curves in Euclidean spaces of arbitrary dimension is investigated.
First, it is shown that the gap between the lowest point of the generalised spectrum 
and the essential spectrum  is positive whenever the cross-section is circular
and the tube is asymptotically straight, untwisted and non-trivially bent.
Second, a Hardy-type inequality is derived for unbent 
and non-trivially twisted tubes. 
\end{abstract}


\maketitle

\section{Introduction}\label{intro}
The interplay between the geometry and spectrum 
of Riemannian manifolds constitutes a traditional area of mathematics.
Physical motivations range from classical vibrational systems
to modern nanostructure devices in quantum mechanics.
Moreover, the study is intrinsically charming
due to the emotional impacts geometric shapes 
have over a person's perception of the world.
The spectrum of the Laplacian in any compact manifold is purely discrete.
On the other hand, non-compact manifolds typically have 
an essential spectrum and the existence of eigenvalues 
is a non-trivial property. 
For non-complete manifolds, the additional challenge in the game 
is the role of boundary conditions.  

The present paper is motivated by an extensive study of 
the Dirichlet Laplacian in a special case of 
non-compact non-complete manifolds: tubes.  
We restrict to the simplest non-trivial situation 
of tubular neighbourhoods of unbounded curves embedded in Euclidean spaces.
Here a strong physical motivation comes from quantum mechanics,
where the Laplacian models the Hamiltonian 
of quantum waveguides \cite{Exner-Kovarik}.

\subsubsection*{Bending is attractive}
Any straight tube has a purely essential spectrum.
Bending it, however, leads to the existence of eigenvalues 
below the essential spectrum of the Dirichlet Laplacian.
This astonishing observation goes back to the pioneering
paper of Exner and \v{S}eba in two dimensions from 1989 \cite{ES}. 
Among the multitude of subsequent results, 
let us highlight the milestones of the generalisation
to three-dimensional tubes via a robust variational proof \cite{GJ,DE},
arbitrary dimensions \cite{cdfk}
and optimal regularity hypotheses \cite{KSed}.
Relying on the quantum-mechanical motivation,
the spectral result can be illustratively interpreted as that an electron
in a curved quantum waveguide gets trapped.

\subsubsection*{Twisting is repulsive}
Finding a way how to geometrically eliminate the discrete eigenvalues 
inspired Ekholm, Kova\v{r}\'ik and one of the present authors 
to establish a Hardy inequality 
for the Dirichlet Laplacian
in twisted three-dimensional tubes \cite{EKK}.
A more robust technique to derive the geometrically induced
Hardy inequalities was later developed in \cite{K6-with-erratum,KZ1}.
Roughly, twisting the waveguide stabilises quantum transport.
An alternative application to the heat flow 
(including the Brownian motion)
can be found in \cite{KZ1,Grillo-Kovarik-Pinchover_2014}. 

\subsection{The nonlinear setting}
The objective of this paper is to generalise 
the spectral-geometric properties of tubes
to the nonlinear setting of the $p$-Laplacian $-\Delta_p$
with $p \in (1,\infty)$
formally acting as 
\begin{equation}\label{operator}
  -\Delta_p u := -\divergence (|\nabla u|^{p-2} \nabla u)
  \,.
\end{equation}  
The linear case mentioned above corresponds to $p=2$, 
while the purely geometric 
setting $p=1$ (the Cheeger problem) is analysed in 
\cite{KP,Leonardi-Pratelli_2016,KLV}.
Our motivation is not only the mathematical curiosity
about the robustness of the effects of bending and twisting, 
but also the relevance of the nonlinear Laplacian in various areas 
of physics and other natural sciences~\cite{Benedikt}.
Moreover, the generalisation is challenging because of
the absence of the powerful spectral theory of
self-adjoint operators whenever $p \not= 2$. 
 
Our ``spectral analysis'' of~\eqref{operator}
is notably inspired by the criticality approach of Pinchover et al.\
\cite{Pinchover_2007,Pinchover-Tintarev_2009,
Pinchover-Regev_2015,Devyver-Pinchover_2016,
Lamberti-Pinchover_2019,Ganguly-Pinchover_2020,dpv}.
Given any open set $\Omega \subset \R^d$ of dimension $d \geq 1$,
we introduce the \emph{spectral threshold}
and the \emph{essential spectral threshold}
by the variational formulae
\begin{equation}\label{ray}
  \lambda_1(\Omega)
  := \inf_{\stackrel[u \not= 0]{}{u\in W_0^{1,p}(\Omega)}} 
  \frac{\displaystyle \int_\Omega |\nabla u(x)|^p \, dx}
  {\displaystyle \int_\Omega |u(x)|^p \, dx}
  \qquad \mbox{and} \qquad
  \lambda_\infty(\Omega) 
  := \sup_{K\Subset \Omega} \lambda_1(\Omega \setminus K)
  \,.
\end{equation}
If~$\Omega$ is bounded, then the infimum is achieved 
and the Euler--Lagrange equation associated with 
the minimisation formula is the quasilinear eigenvalue problem
\begin{equation}\label{quasi}
  -\Delta_p u = \lambda_1(\Omega) |u|^{p-2} u
\end{equation}
in~$\Omega$,
subject to Dirichlet boundary conditions $u=0$ on~$\partial\Omega$;  
moreover, $\lambda_\infty(\Omega) = \infty$.
In general, $\lambda_1(\Omega)$ and $\lambda_\infty(\Omega)$
extend to $p \not= 2$ the well-known variational characterisations of 
the lowest point in the spectrum (Rayleigh--Ritz)
and the essential spectrum (Persson), respectively,
of the self-adjoint Dirichlet Laplacian $-\Delta_2$ in $L^2(\Omega)$.

Whenever, the infimum in~\eqref{ray} is achieved by 
a function $u\in W_0^{1,p}(\Omega)$, 
we call it the \emph{first eigenfunction} 
(or \emph{ground state}) of~$\Omega$.
Then the spectral threshold $\lambda_1(\Omega)$ is also called 
the \emph{first} (or \emph{principal}) \emph{eigenvalue} of~$\Omega$.
For unbounded~$\Omega$, the existence of the ground state
is a highly non-trivial property.
In particular, the positivity of the \emph{essential spectral gap}
$\lambda_\infty(\Omega) - \lambda_1(\Omega)$ 
generalises the existence of (discrete) eigenvalues
below the essential spectrum of $-\Delta_2$.

The class of domains we are interested in in this paper
are deformations of the \emph{straight tube} $\Omega_0 := \R \times \omega$,
where $\omega \subset \R^{d-1}$ with $d \geq 2$ 
is an arbitrary bounded open connected set.  
It is easy to see that  
\begin{equation}\label{straight} 
  \lambda_1(\Omega_0)=\lambda_1(\omega)=\lambda_\infty(\Omega_0)
  \,,
\end{equation}
so the essential spectral gap is zero in this case.
Our goal is to analyse the influence of 
bending and twisting of~$\Omega_0$
on the spectral threshold and the essential spectral threshold. 

\subsection{The geometric framework}
Let $\Gamma: \mathbb R \to \mathbb R^d$ be a $C^{1,1}$-smooth unit-speed curve. 
Then $T:=\Gamma'$ is a unit tangent vector field along~$\Gamma$
and $\kappa:=|\Gamma''|$ is the (locally bounded) curvature function of $\Gamma$.
There exist (almost everywhere differentiable) 
unit normal vector fields $N_1,\dots, N_{d-1}:\R\to\R^d$ such that
\begin{equation}\label{kij}
\begin{pmatrix}
T\\N_1\\\vdots\\ N_{d-1}\end{pmatrix}'=\begin{pmatrix} 
0&\kappa_1&\dots&\kappa_{d-1}\\
-\kappa_1&0&\dots&0\\
\vdots&\vdots&& \vdots\\
-\kappa_{d-1}&0&\dots&0
\end{pmatrix}
\begin{pmatrix}
T\\N_1\\\vdots\\ N_{d-1}\end{pmatrix}
\,,
\end{equation}
where $\kappa_1,\dots,\kappa_{d-1}:\R\to\R$ are (locally bounded) functions 
such that $\kappa_1^2+\dots+\kappa_{d-1}^2=\kappa^2$.
Let us consider a one-parametric
family of rotation matrices
$$
  R: \R \to \mathsf{SO}(d-1)
  \,,
$$
which we assume to be differentiable with 
$R' \in L_\mathrm{loc}^\infty(\R;\R^{(d-1)\times(d-1)})$.
Rotating the normal vector fields $N_1,\dots,N_{d-1}$ via~$R$,
we obtain an arbitrary frame
$(T,R_{1\mu}N_\mu,\dots,R_{d-1\mu}N_\mu)$ of~$\Gamma$.
Here the Einstein summation convention is adopted,
with the range of Greek indices being $1,\dots,d-1$.
Then a general \emph{bent twisted tube} about~$\Gamma$ is obtained by 
\begin{equation}\label{tube}
  \Omega_{\kappa,R} := \big\{
  \Gamma(s) + t_\mu \, R_{\mu\nu}(s) \, N_\nu 
  : \ (s,t) \in \R \times \omega 
  \big\}
  \,.
\end{equation}
Our standing hypothesis is that~$\Omega_{\kappa,R}$ does not overlap itself,
which particularly involves the necessary condition that 
$\kappa \in L^\infty(\R)$ and 
\begin{equation}\label{Ass}
  a \, \|\kappa\|_{L^\infty(\R)}  < 1  
\end{equation}
with $\displaystyle a := \sup_{t \in \omega} |t|$. 

We say that $\Omega_{\kappa,R}$ is \emph{unbent} or \emph{untwisted}
if $\kappa=0$ or $R'=0$, respectively. 
The former implies that~$\Gamma$ is a straight line,
while the latter means that 
the cross-section~$\omega$
is translated along~$\Gamma$ with respect to 
a \emph{relatively parallel} frame 
(i.e., the normal components rotate along~$\Gamma$ 
only whatever amount is necessary to remain normal,
namely, their derivative stays tangential).
An unbent untwisted tube 
is straight in the sense that it is congruent to
$\Omega_0 = \Omega_{0,I}$,
where~$I$ is the identity matrix.

The terminology is not perfect because straight tubes are considered
as a particular situation of bent twisted tubes.
For this reason, we say that $\Omega_{\kappa,R}$
is \emph{non-trivially bent} if $\kappa\not=0$.
However, even if $R'\not=0$,
it might happen that $\Omega_{0,R}$ is congruent to~$\Omega_0$.
In fact, this is the case whenever the cross-section~$\omega$ 
is \emph{circular}, i.e.,
$\omega$~is
the ball $B_a(0)$ of radius $a>0$ centred at the origin of~$\R^{d-1}$
or $\omega = B_a(0) \setminus \overline{B_{a_0}(0)}$ is a spherical shell
of radii $0 < a_0 < a$ centred at the origin.
Therefore, the property of a tube being \emph{non-trivially twisted} 
requires an extra hypothesis about the asymmetry 
of the cross-section~$\omega$ with respect to the rotations~$R$.
A discovery of this paper is that, in all dimensions, 
the right definition reads
\begin{equation}\label{twist}
  f_\mu\partial_{t_\mu} \phi_1 \not= 0
\end{equation}
as an identity in~$\Omega_0$.  
Here~$\phi_1$ is the first eigenfunction of~$\omega$ 
and $f_\mu(s,t) := t_\alpha R_{\alpha\beta}'(s) R_{\mu\beta}(s)$. 

\subsection{The main results}
Our first result is about \emph{asymptotically straight} tubes 
characterised by the vanishing of bending and twisting at infinity:
\begin{equation}\label{Ass.decay}
  \lim_{|s|\to\infty} \kappa(s) = 0
  \qquad \mbox{and} \qquad
  \lim_{|s|\to\infty} 
  \int_\omega |f_\mu(s,t) \partial_{t_\mu} \phi_1(t)|^p \, dt 
  = 0
  \,.
\end{equation}
Of course, a sufficient condition to ensure the validity
of the second limit is that $R'(s) \to 0$ as $|s| \to \infty$
(in any matrix topology).
\begin{thm}[stability of~$\lambda_\infty$]\label{Thm.ess}
If~\eqref{Ass.decay} holds, then
$
  \lambda_\infty(\Omega_{\kappa,R}) = \lambda_1(\omega) 
$.  
\end{thm}
\noindent
In particular, the second equality of~\eqref{straight}
follows as a very special case.

Our second result is that the spectral threshold diminishes
whenever the tube is non-trivially bent and untwisted.
Unfortunately, we are able to prove it only if 
the cross-section~$\omega$ is circular.  
Recall that $\Omega_{\kappa,R}=\Omega_{\kappa,I}$ in this case,
even if $R'\not=0$.
\begin{thm}[bending]\label{Thm.disc}
If $\kappa \not= 0$ and $\omega$ is circular, then
$
  \lambda_1(\Omega_{\kappa,I}) < \lambda_1(\omega) 
$.  
\end{thm}

As a consequence of Theorems~\ref{Thm.ess} and~\eqref{Thm.disc},
we get the ultimate result about the positivity of the essential
spectral gap.
\begin{cor}\label{Corol.disc}
If $\kappa \not= 0$, 
$\omega$ is circular and~\eqref{Ass.decay} holds, then
$
  \lambda_1(\Omega_{\kappa,I})
  < \lambda_\infty(\Omega_{\kappa,I}) 
$.  
\end{cor}
\noindent
We leave as an open problem (see Remark~\ref{Rem.circular})
whether the result holds for arbitrary cross-sections
(unless $p=2$ when the general validity is well known).

Finally, for unbent non-trivially twisted tubes 
we establish a Hardy inequality.
\begin{thm}[twisting]\label{Thm.Hardy}
If $\kappa=0$ and~\eqref{twist} holds, 
then there exists a positive 
continuous function $\rho: \Omega_{0,R} \to \R$ such that 
\begin{equation}\label{Hardy.intro}
  \forall u\in W_0^{1,p}(\Omega_{0,R})
  \,, \qquad
  \int_{\Omega_{0,R}} |\nabla u|^p \, dx
  - \lambda_1(\omega) \int_{\Omega_{0,R}} |u|^p \, dx
  \geq \int_{\Omega_{0,R}} \rho \, |u|^p \, dx
  \,.
\end{equation} 
\end{thm}

Note that the theorem is void if $d=2$,
because there is no twisting for a two-dimensional strip.
If $d = 3$ and $p=2$, the existence of a Hardy inequality 
is known from \cite{EKK,K6-with-erratum,KZ1},
however, a positive weight~$\rho$ was established only 
for compactly supported~$R'$. 
Here we provide a robust existence of the Hardy inequality
under the minimal hypothesis. 
What is more, we prove the Hardy inequality in all dimensions  
$d \geq 3$ and $p \in (1,\infty)$.  
The result is completely new in higher dimensions $d \geq 4$ 
even in the linear case $p=2$.

The result Theorem~\ref{Thm.Hardy} is highly non-trivial because 
there is no Hardy inequality in straight tubes~$\Omega_0$.
Indeed, the shifted operator 
$-\Delta_p - \lambda_1(\omega)$ in~$\Omega_0$ 
is \emph{critical} in the sense that 
the spectral threshold of $-\Delta_p - \lambda_1(\omega) + V$ in~$\Omega_0$
is negative whenever the perturbation $V \in C_0^\infty(\Omega_0)$ 
is non-positive and non-trivial (yet arbitrarily small),
see Proposition~\ref{Prop.critical}.
Theorem~\ref{Thm.Hardy} says that 
$-\Delta_p - \lambda_1(\omega)$ in~$\Omega_{0,R}$ is \emph{subcritical} 
whenever~$\Omega_{0,R}$ is non-trivially twisted, 
in the sense of the stability
that the spectral threshold of $-\Delta_p - \lambda_1(\omega) + V$ 
in~$\Omega_{0,R}$ remains zero whenever $V \in C_0^\infty(\Omega_{0,R})$ is small.
For non-trivially bent untwisted tubes,
the operator $-\Delta_p - \lambda_1(\omega)$ in~$\Omega_{\kappa,I}$
may be understood as  
\emph{supercritical} under the hypotheses of Theorem~\ref{Thm.disc},
because $\lambda_1(\Omega_{\kappa,I}) -\lambda_1(\omega)$ is negative
even if $V=0$.

\subsection{Possible extensions} 
As an interesting direction of possible future research,
let us mention the replacement of the Dirichlet boundary conditions
by \emph{Robin} boundary conditions in the spirit of 
\cite{Kovarik-Pankrashkin_2017}.
Apart from two-dimensional studies 
\cite{DKriz2,FK3,K5,Novak_2016,Rossini_2019}
or thin-width asymptotics  
\cite{BMT2,Oliveira-Rossini_2022},
a detailed spectral-geometric analysis of higher-dimensional Robin waveguides
remains open even in the linear case $p=2$.

Another interesting extension would be 
to add \emph{magnetic}   fields to the $p$-Laplacian
\cite{KR,CKLL}. 

The present paper is essentially concerned with the quasi-linear
eigenvalue problem~\eqref{quasi} 
(to be properly understood through~\eqref{ray}).
This is a special (so-called \emph{$p$-linear}) 
case of the more general problem
\begin{equation*} 
  -\Delta_p u = \lambda \, |u|^{q-2} u
  \,.
\end{equation*}
The other cases $q>p$ and $q < p$ are known as 
\emph{$p$-sublinear} and \emph{$p$-superlinear}, respectively.
There is an extensive literature 
on existence of solutions both in the variational and
non-variational frameworks, including more general operators, 
see, e.g., \cite{AP, blsch, bjde, GPeral, ruiz} and references therein.
We consider an extension of the present spectral-geometric study 
to the more general types of nonlinearities as 
yet another interesting future project.

\subsection{The organisation} 
The paper is structured as follows. 
In Section~\ref{prel} we comment more on the geometric setting
and implement the usual passage to the curvilinear ``coordinates'' $(s,t)$,
which is the main strategy to deal with curved quantum waveguides.
Straight tubes are considered in Section~\ref{stra2},
where we establish the first equality of~\eqref{straight}
as well as Proposition~\ref{Prop.critical}
about the criticality of~$\Omega_0$.
Theorems~\ref{Thm.ess}, \ref{Thm.disc} and~\ref{Thm.Hardy} 
are established in Sections~\ref{asy3}, \ref{bent3} and~\ref{Sec.twist}, respectively.
In the last section concerned with twisted tubes,
we also comment more on hypothesis~\eqref{twist}.

\section{Preliminaries}\label{prel}
We refer to \cite{KSed,kz} on the geometry of curves
under the present minimal hypotheses. 
The unit-speed hypothesis about~$\Gamma$ means that 
$|\Gamma'(s)| = 1$ for every $s \in \R$,
which can be always achieved by a suitable (arc-length) change of parameterisation.
Note that the \emph{relatively parallel adapted frame}
$(T, N_1,\dots, N_{d-1})$ 
differs from the customarily used \emph{Frenet frame}.
The latter requires an extra (classically $C^d$-smoothness) regularity of~$\Gamma$,
and moreover, curves with vanishing curvature somewhere must be excluded.
The relatively parallel adapted frame is uniquely defined
modulo the choice of initial conditions 
$N_j(s_0) = N_j^0(s_0) \in \R^d$ for some $s_0 \in \R$, 
which also fixes the curvature functions~$\kappa_j$,
with $j \in \{1,\dots,d-1\}$.
The shape of a bent untwisted tube~$\Omega_{\kappa,I}$ 
therefore depends on the choice of the relatively parallel adapted frame,
unless the cross-section~$\omega$ is circular. 

For the usual passage to the curvilinear ``coordinates'' $(s,t)$
when dealing with curved tubes, 
we refer to the geometrically oriented references \cite{cdfk,kz,K-Padova}. 
Let us introduce the mapping 
$\mathscr{L} : \R\times \R^{d-1}\to \R^d$ defined by 
$$  
  \mathscr{L}(s, t) := \Gamma(s) + t_\mu \, R_{\mu\nu}(s) \, N_\nu 
  \,,
$$
so that $\Omega_{\kappa,R} = \mathscr{L}(\Omega_0)$.
It is convenient to think of~$\Omega_{\kappa,R}$ 
as the Riemannian manifold~$\Omega_0$
equipped with the induced metric 
$
  g := (\nabla \mathscr{L}) \cdot (\nabla \mathscr{L})^T
$,
where the dot denotes the matrix multiplication in~$\R^d$. 
Using~\eqref{kij} and the orthogonality of~$R$, 
it is straightforward to check that the metric reads 
$$
  g =
  \begin{pmatrix}
    f^2 + f_\mu f_\mu & f_1 & f_2 & \dots & f_{d-1} \\
    f_1 & 1   & 0   & \dots & 0 \\
    f_2 & 0   & 1   & \dots & 0 \\
    \vdots & \vdots & \vdots & \ddots & \vdots \\
    f_{d-1} & 0 & 0 & \dots & 1
  \end{pmatrix}
  , \qquad  
  \det(g) := f^2 
  \,,
$$
where
$$
  \begin{aligned}
    f(s,t) :=& \ 1-t_\alpha R_{\alpha\beta}(s)\kappa_\beta(s) \,,    
    \\
    f_\mu(s,t) :=& \ t_\alpha R_{\alpha\beta}'(s) R_{\mu\beta}(s) \,.
  \end{aligned}
$$ 
Note that 
$
  (f_\mu f_\mu)(s,t) 
  = t_\alpha R_{\alpha\beta}'(s) t_\nu R_{\nu\beta}'(s) 
$
by the orthogonality of~$R$. 
From the basic hypothesis~\eqref{Ass}, 
it follows that the Jacobian of~$\mathscr{L}$ satisfies
\begin{equation}\label{est}
0<1-a\, \|\kappa\|_{L^\infty(\R)}
\leq f(s,t) \leq
1+a\, \|\kappa\|_{L^\infty(\R)} < \infty
\end{equation}
for every $(s,t)\in\Omega_0$.
Consequently, $(\Omega_0,g)$ is a Riemannian manifold
provided that~\eqref{Ass} holds.
More specifically, 
$\mathscr{L}:\Omega_0\to\Omega_{\kappa,R}$ is a local $C^{0,1}$-diffeomorphism
under the assumption~\eqref{Ass}
(cf.~\cite[Prop.~2.2]{kz}).
To make it a global diffeomorphism, 
one needs to additionally assume that~$\mathscr{L}$ is injective.
The inverse metric is given by
$$
  g^{-1} =
  \begin{pmatrix}
    1 & -f_1 & -f_2 & \dots & -f_{d-1} \\
    -f_1 & 1 + f_1^2   & f_1 f_2   & \dots & f_1 f_{d-1} \\
    -f_2 & f_2 f_1   & 1 + f_2^2   & \dots & 0 \\
    \vdots & \vdots & \vdots & \ddots & \vdots \\
    -f_{d-1} & f_{d-1} f_1 & f_{d-1} f_2 & \dots & 1 + f_{d-1}^2
  \end{pmatrix}
  .
$$

Passing to the curvilinear coordinates
in the integrals 
in the variational characterisation~\eqref{ray} 
by the change of trial function $\psi := u \circ \mathscr{L}$, 
it is straightforward to verify that 
\begin{equation}\label{ray.bis} 
  \lambda_1(\Omega_{\kappa,R})
  = \inf_{\stackrel[\psi \not= 0]{}{\psi\in W_0^{1,p}(\Omega_0,g)}} 
  \frac{\displaystyle Q[\psi]}
  {\displaystyle \|\psi\|^p}
  =: \lambda_1(\Omega_0,g)
  \,,
\end{equation}
where
\begin{equation}\label{form} 
\begin{aligned}
  Q[\psi] &:= \int_{\Omega_0} 
  \left(
  \left|\frac{(\partial_s-f_\mu(s,t)\partial_{t_\mu})\psi(s,t)}{f(s,t)}\right|^2 
  + |\nabla_{\!t} \psi(s,t)|^2
  \right)^{p/2} f(s,t) \, ds \, dt  
  \,, \\
  \|\psi\| & :=
  \left(
  \int_{\Omega_0} |\psi(s,t)|^p \, f(s,t) \, ds \, dt 
  \right)^{1/p}
  \,,
\end{aligned}  
\end{equation}
and $W_0^{1,p}(\Omega_0,g)$ denotes the closure of $C_0^\infty(\Omega_0)$
with respect to the norm $(Q[\psi] + \|\psi\|^p)^{1/p}$.
By virtue of~\eqref{est}, it is straightforward to verify that 
$W_0^{1,p}(\Omega_0,g) = W_0^{1,p}(\Omega_0)$
provided that $\kappa \in L^\infty(\R)$ (always assumed)
and $R' \in L^\infty(\R;\R^{(d-1)\times(d-1)})$
(not necessarily assumed).

\begin{rem}
It is evident that~\eqref{ray.bis} is well defined merely
under the hypothesis~\eqref{Ass}.
Therefore, if we take~\eqref{ray.bis} 
as the very definition of $\lambda_1(\Omega_{\kappa,R})$ 
and abandon the interpretation of~$\Omega_{\kappa,R}$ as a non-self-intersecting tube,
all the results in this paper hold without the extra assumption
that~$\mathscr{L}$ is injective.
In this more general approach, the essential spectral threshold 
should be interpreted as 
\begin{equation}\label{ray.bis.ess} 
  \lambda_\infty(\Omega_{\kappa,R})
  = \sup_{K \Subset \Omega_0} \lambda_1(\Omega_0\setminus K,g)
  \,.
\end{equation}
\end{rem}

To handle $Q[\psi]$, the following elementary observation will be widely used.
\begin{lem}\label{Lem.crucial}
For any non-negative numbers $a,b,q$, one has 
\begin{equation}\label{crucial}
  (a+b)^q \leq \alpha^q a^q + \beta^q b^q
  \,,
\end{equation}
where $\alpha,\beta$ are any positive numbers satisfying 
$\frac{1}{\alpha} + \frac{1}{\beta} = 1$.
\end{lem}
\begin{proof}
If $q \leq 1$, one has the better inequality 
$(a+b)^q \leq a^q + b^q$, from which~\eqref{crucial} follows
by the fact that necessarily $\alpha,\beta > 1$.
In any case, the claim is achieved by arguing that 
either $a+b \leq \alpha a$ or $a+b \leq \beta b$ holds.
This is true since otherwise one would get the contradiction that
$(\frac{1}{\alpha} + \frac{1}{\beta})(a+b) > a+b$.
\end{proof} 

\section{Straight tubes}\label{stra2}
If $\kappa=0$ and $R'=0$, then $\Gamma$~is a straight line
and the relatively parallel adapted frame is actually parallel
(i.e., constant along~$\Gamma$).
Consequently, $\Omega_{\kappa,R}$ coincides with 
$\Omega_0=\mathbb R\times \omega$ up to congruence. 
Our goal is to establish the first equality of~\eqref{straight} 
for any cross-section~$\omega$. 
Since any straight tube is necessarily asymptotically straight,
the proof of the second equality of~\eqref{straight} 
is postponed to Section~\ref{asy3}. 
 
First of all, since~$\omega$ is assumed to be bounded and connected,
it is well known \cite{lind,Belloni-Kawohl_2002,Kawohl-Lindqvist}
that $\lambda_1(\omega)$ is a simple eigenvalue of 
the Dirichlet $p$-Laplacian in~$\omega$.
More specifically, there exists a unique 
(up to a constant multiple) 
positive minimiser $\phi_1 \in W_0^{1,p}(\omega)$ 
of the minimisation problem in~\eqref{ray}
(with~$\Omega$ being replaced by~$\omega$).
We choose it normalise to~$1$ in $L^p(\omega)$,
i.e., $\int_\omega |\phi_1(t)|^p \, dt = 1$.
The variational characterisation of $\lambda_1(\omega)$ 
yields the Poincar\'e inequality
\begin{equation}\label{Poincare} 
  \forall \phi \in W_0^{1,p}(\omega) 
  \,, \qquad
  \int_\omega |\nabla\phi(t)|^p \, dt
  \geq \lambda_1(\omega) \int_\omega |\phi(t)|^p \, dt
  \,.
\end{equation}
Note that $\lambda_1(\omega)$ is necessarily positive.
We assume no regularity hypotheses about~$\omega$,
unless otherwise stated.

\begin{Prop}\label{Prop.straight} 
One has 
$$
  \lambda_1(\Omega_0) = \lambda_1(\omega)
  \,.
$$ 
\end{Prop}
\begin{proof}
If $\kappa=0$ and~$R'=0$, then $f=1$ and $f_\mu=0$. 
Consequently,
$$
\begin{aligned}
  Q[\psi] 
  &=  \int_{\Omega_0} 
  \left(
  |\partial_s\psi|^2 
  + |\nabla_{\!t} \psi|^2
  \right)^{p/2} ds \, dt 
  \\
  &\geq 
  \int_{\Omega_0} 
  |\nabla_{\!t} \psi|^p
  \, ds \, dt 
  \\
  &\geq \lambda_1(\omega)
   \int_{\Omega_0} 
  |\psi|^p
  \, ds \, dt 
  = \lambda_1(\omega) \, \|\psi\|^p
  \,,
\end{aligned}   
$$
where the last inequality follows from~\eqref{Poincare} 
and the Fubini theorem.
(For simplicity, we suppress the arguments of functions 
in the integrals from now on.)

To prove the opposite inequality, 
it is enough to construct a sequence 
$(\psi_n)_{n=1}^\infty \subset W_0^{1,p}(\Omega_0)$
such that 
\begin{equation}\label{Weyl} 
  R[\psi_n] := 
  \frac{Q[\psi_n] - \lambda_1(\omega) \, \|\psi_n\|^p}{\|\psi_n\|^p}
  \xrightarrow[n \to \infty]{} 0
  \,.
\end{equation}
To this purpose, let $(\varphi_n)_{n=1}^\infty \subset W_0^{1,p}(\R)$
be defined by 
\begin{equation}\label{one}
\begin{aligned}
  \varphi_n(s) :=
  \begin{cases}
    1 & \mbox{if} \quad |s| \leq n \,, \\
    \displaystyle
    2-\frac{|s|}{n} & \mbox{if} \quad |s| \in (n,2n) \,, \\ 
    0 & \mbox{if} \quad |s| \geq 2n \,. \\ 
  \end{cases}
\end{aligned}
\end{equation}
Note that $\varphi_n \to 1$ pointwise as $n \to \infty$ and 
\begin{equation}\label{fnp}
  \int_{\mathbb R} |\varphi_n'(s)|^\xi \, ds
  = \frac{2}{n^{\xi-1}} 
  \xrightarrow[n \to \infty]{} 0
\end{equation}
whenever $\xi > 1$. 
Set $\psi_n(s,t) := \varphi_n(s) \phi_1(t)$.
By Lemma~\ref{Lem.crucial},
$$
\begin{aligned}
  Q[\psi_n] 
  &\leq  
  \alpha^{p/2} \int_{\Omega_0} 
  |\partial_s\psi_n|^p \, ds \, dt
  +  \beta^{p/2} \int_{\Omega_0} 
  |\nabla_{\!t} \psi_n|^p \, ds \, dt 
  \\
  &= \alpha^{p/2} \int_{\R} 
  |\varphi_n'|^p \, ds
  +  \beta^{p/2} \, \lambda_1(\omega) \int_{\R} 
  |\varphi_n|^p \,
  ds  
  \,,
\end{aligned}  
$$
where the equality follows by the fact  
that equality holds in~\eqref{Poincare} 
if (and only if) $\phi = \phi_1$ 
and by the normalisation of~$\phi_1$.
At the same time, 
$
  \|\psi_n\|^p = \int_{\R} |\varphi_n|^p \,  ds  
$.
Consequently,
$$
  R[\psi_n] \leq 
  \alpha^{p/2} \, 
  \frac{\displaystyle\int_{\R} |\varphi_n'|^p \, ds}
  {\displaystyle\int_{\R} |\varphi_n|^p \, ds}
  + (\beta^{p/2} -1) \, \lambda_1(\omega) 
  \,.
$$
By~\eqref{fnp}, it follows that 
$$
  \lim_{n \to \infty} R[\psi_n] 
  \leq (\beta^{p/2}-1) \, \lambda_1(\omega) 
  \,,
$$
where $\beta > 1$ can be made arbitrarily close to~$1$.
\end{proof}

Straight tubes are critical in the sense of 
the following instability of $-\Delta_p$
with respect to small perturbations.
\begin{Prop}\label{Prop.critical}
Let $V \in C_0^\infty(\Omega_0)$ be non-positive and non-trivial.
Then
$$
\begin{aligned}
  \lambda_1^V(\Omega_0) := & \
  \inf_{\stackrel[\psi \not= 0]{}{\psi\in W_0^{1,p}(\Omega_0)}} 
  \frac{\displaystyle  \int_{\Omega_0} |\nabla\psi|^p \, ds \, dt
  + \int_{\Omega_0} V |\psi|^p \, ds \, dt}
  {\displaystyle \|\psi\|^p}
  \\
  < & \ 
  \inf_{\stackrel[\psi \not= 0]{}{\psi\in W_0^{1,p}(\Omega_0)}} 
  \frac{\displaystyle  \int_{\Omega_0} |\nabla\psi|^p \, ds \, dt}
  {\displaystyle \|\psi\|^p}
  = \lambda_1(\Omega_0) = \lambda_1(\omega) 
\end{aligned} 
$$  
\end{Prop}
\begin{proof}
The last but one inequality is the definition of $\lambda_1(\Omega_0)$,
while the last inequality is~\eqref{straight} (see Proposition~\ref{Prop.straight}). 
The main claim is the strict inequality.
To prove it, it is enough to find a (trial) function
$\psi \in  W_0^{1,p}(\Omega_0)$ for which
$$
  Q_1^V[\psi] := \int_{\Omega_0} |\nabla\psi|^p \, ds \, dt
  + \int_{\Omega_0} V \, |\psi|^p \, ds \, dt
  - \lambda_1(\omega) \int_{\Omega_0} |\psi|^p \, ds \, dt
  < 0 \,.
$$ 
It is achieved by the trial function $\psi_n(s,t) := \varphi_n(s) \phi_1(t)$,
where the sequence~$\varphi_n$ is defined in~\eqref{one}. 
Indeed, if $p \leq 2$, then
$$
\begin{aligned}
  Q_1^V[\psi_n] 
  &\leq \int_{\Omega_0} 
  |\partial_s\psi_n|^p
  \, \, ds \, dt 
  + \int_{\Omega_0} V \, |\psi_n|^p \, ds \, dt
  \\
  &= \int_{\Omega_0} 
  |\varphi_n'|^p
  \, \, ds 
  + \int_{\Omega_0} V \, |\varphi_n|^p \, |\phi_1|^p \, ds \, dt
  \xrightarrow[n \to \infty]{} 
  \int_{\Omega_0} V(s,t) \, |\phi_1(t)|^p \, ds \, dt 
  < 0
  \,,
\end{aligned}  
$$ 
where we have used  
that equality holds in~\eqref{Poincare} 
if (and only if) $\phi = \phi_1$  
and the normalisation of~$\phi_1$.
The convergence holds due to~\eqref{fnp}
and since~$\varphi_n$ converges to~$1$ pointwise. 
If $p>2$, Lemma~\ref{Lem.crucial} yields
$$
\begin{aligned}
  Q_1^V[\psi_n] 
  &\leq  
  \alpha^{p/2} \int_{\Omega_0} 
  |\partial_s\psi_n|^p  \, ds \, dt
  +  (\beta^{p/2}-1) \, \lambda_1(\omega) \int_{\Omega_0} 
  |\psi_n|^p \, ds \, dt 
  + \int_{\Omega_0} V \, |\psi_n|^p \, ds \, dt
  \\
  &= \alpha^{p/2} 
  \int_{\R} |\varphi_n'|^p \, ds 
  +  (\beta^{p/2}-1) \, \lambda_1(\omega) \int_{\R} |\varphi_n|^p \, ds 
  + \int_{\Omega_0} V \, |\varphi_n|^p \, |\phi_1|^p \, ds \, dt
  \\
  &\leq  2\alpha^{p/2} \,
  \frac{1}{n^{p-1}} 
  +  (\beta^{p/2}-1) \, \lambda_1(\omega) \, 4 n
  + \int_{\Omega_0} V \, |\varphi_n|^p \, |\phi_1|^p \, ds \, dt
  \,,
\end{aligned} 
$$
where the last inequality employs~\eqref{fnp} 
and the fact that the function~$\varphi_n$ satisfies
$0 \leq \varphi_n \leq 1$ on $(-2n,2n)$
and that it is zero elsewhere.
Choosing $n$-dependent~$\beta$ 
(and thus also $n$-dependent conjugate~$\alpha$), 
for instance ($n \geq 2$), 
\begin{equation}\label{depend}
  \beta := 1 + (n \log n)^{-1}
  \qquad \mbox{(which implies $\alpha = 1 + n \log n$)}
\end{equation}
it is straightforward to check that 
$$
   \lim_{n \to \infty} Q_1^V[\psi_n] 
   \leq \int_{\Omega_0} V(s,t) \, |\phi_1(t)|^p \, ds \, dt 
   < 0
$$ 
in this case as well.
In summary, there exists a natural number~$n_0$
such that $Q_1^V[\psi_n]$ is negative for all $n \geq n_0$.
\end{proof}

\section{Asymptotically straight tubes}\label{asy3}
In this section we establish Theorem~\ref{Thm.ess}.
To simplify the presentation, we divide the proof into two steps.
First, we show that the essential spectral threshold 
satisfies the required lower bound
if the tube is merely asymptotically unbent.

\begin{Prop}\label{Prop.lower}
If $\displaystyle \lim_{|s|\to\infty} \kappa(s) = 0$, then
$$
  \lambda_\infty(\Omega_{\kappa,R}) \geq \lambda_1(\omega) 
  \,.
$$  
\end{Prop}
\begin{proof}
By the definition of $\lambda_\infty(\Omega_{\kappa,R})$ 
given in~\eqref{ray.bis.ess}, 
one has 
$
  \lambda_\infty(\Omega_{\kappa,R}) 
  \geq \lambda_1(\Omega_0 \setminus K,g)
$ 
for any ``trial'' compact subset~$K$ of~$\Omega_0$.  
For any positive numbers~$\eps$ (small) and~$l$ (large), we set 
$$
  K := [-l,l] \times \overline{\omega_\eps}
  \qquad \mbox{with} \qquad
  \omega_\eps := \{t \in \omega : \dist(t,\partial\omega) > \eps \}
  \,.
$$
Let $\psi \in W_0^{1,p}(\Omega_0 \setminus K)$ be arbitrary.
Neglecting the ``longitudinal energy'' of $Q[\psi]$ in~\eqref{form}
and recalling~\eqref{est}, 
one has 
$$
\begin{aligned}
  Q[\psi] &\geq 
  \int_{\Omega_0} |\nabla_{\!t} \psi|^p \, f \, ds \, dt 
  \\
  &= \int_{[-l,l]} \int_{\omega\setminus\overline{\omega_\eps}} 
  |\nabla_{\!t} \psi|^p \, f \, dt \, ds 
  + \int_{\R\setminus[-l,l]} \int_{\omega} 
  |\nabla_{\!t} \psi|^p \, f \, dt \, ds 
  \\
  &\geq (1-a\,\|\kappa\|_{L^\infty(\R)}) 
  \int_{[-l,l]} \int_{\omega\setminus\overline{\omega_\eps}} 
  |\nabla_{\!t} \psi|^p \, dt \, ds 
  \\
  & \quad + (1-a\,\|\kappa\|_{L^\infty(\R\setminus[-l,l])}) 
  \int_{\R\setminus[-l,l]} \int_{\omega} 
  |\nabla_{\!t} \psi|^p \, dt \, ds 
  \\
  &\geq (1-a\,\|\kappa\|_{L^\infty(\R)}) 
  \, \lambda_1(\omega\setminus\overline{\omega_\eps})
  \int_{[-l,l]} \int_{\omega\setminus\overline{\omega_\eps}} 
  |\psi|^p \, dt \, ds 
  \\
  & \quad + (1-a\,\|\kappa\|_{L^\infty(\R\setminus[-l,l])}) 
  \, \lambda_1(\omega)
  \int_{\R\setminus[-l,l]} \int_{\omega} 
  |\psi|^p \, dt \, ds 
  \\
  &\geq 
  \frac{1-a\,\|\kappa\|_{L^\infty(\R)}}
  {1+a\,\|\kappa\|_{L^\infty(\R)}} 
  \, \lambda_1(\omega\setminus\overline{\omega_\eps})
  \int_{[-l,l]} \int_{\omega\setminus\overline{\omega_\eps}} 
  |\psi|^p \, f \, dt \, ds 
  \\
  & \quad + 
  \frac{1-a\,\|\kappa\|_{L^\infty(\R\setminus[-l,l])}}
  {1+a\,\|\kappa\|_{L^\infty(\R\setminus[-l,l])}} 
  \, \lambda_1(\omega)
  \int_{\R\setminus[-l,l]} \int_{\omega} 
  |\psi|^p \, f \, dt \, ds 
  \\
  &\geq 
  \min\left\{
  \frac{1-a\,\|\kappa\|_{L^\infty(\R)}}
  {1+a\,\|\kappa\|_{L^\infty(\R)}} 
  \, \lambda_1(\omega\setminus\overline{\omega_\eps})
  ,
  \frac{1-a\,\|\kappa\|_{L^\infty(\R\setminus[-l,l])}}
  {1+a\,\|\kappa\|_{L^\infty(\R\setminus[-l,l])}} 
  \, \lambda_1(\omega)
  \right\}
  \|\psi\|^p
  \,.
\end{aligned}  
$$
Note that $\lambda_1(\omega\setminus\overline{\omega_\eps}) \to \infty$
as $\eps \to 0$. 
To see it, one can argue through the Faber--Krahn inequality 
for the Dirichlet $p$-Laplacian
(see, e.g., \cite{Chorwadwala-Mahadevan-Toledo} and references therein)
and the scaling 
$\lambda_1(B_\epsilon(0)) = \epsilon^{-p} \lambda_1(B_1(0))$.
Consequently, the minimum equals the second constant 
for all sufficiently small~$\eps$, so we have established the bound
$$
  \lambda_\infty(\Omega_{\kappa,R}) 
  \geq 
  \frac{1-a\,\|\kappa\|_{L^\infty(\R\setminus[-l,l])}}
  {1+a\,\|\kappa\|_{L^\infty(\R\setminus[-l,l])}} 
  \ \lambda_1(\omega)
  \,.
$$
Since the fraction tends to~$1$ as $l \to \infty$,
the desired lower bound follows.
\end{proof}

Now we turn to the upper bound.
For simplicity, let us denote
$$
  r(s) :=
  \int_\omega |f_\mu(s,t) \partial_{t_\mu} \phi_1(t)|^p \, dt 
  \,.
$$

\begin{Prop}\label{Prop.upper}
If $\displaystyle \lim_{s\to\infty} \kappa(s) = 0$
and $\displaystyle \lim_{s\to\infty} r(s) = 0$, 
then
$$
  \lambda_\infty(\Omega_{\kappa,R}) \leq \lambda_1(\omega) 
  \,.
$$  
\end{Prop}
\begin{proof}
Fix any $K\Subset\Omega_0$ and let us define 
$\psi_n(s,t):=\tilde{\varphi}_n(s) \phi_1(t)$
with $\tilde{\varphi}_n(s) := \varphi_n(s-n^2)$,
where~$\varphi_n$ is the sequence defined in~\eqref{one}. 
As in the proof of Proposition~\ref{Prop.straight},
$\tilde{\varphi}_n \to 1$ pointwise as $n \to \infty$
and~\eqref{fnp} holds for~$\varphi_n$ 
being replaced by~$\tilde{\varphi}_n$ as well.
Moreover, $\tilde{\varphi}_n$ is ``localised at''~$\infty$ meaning that
$\inf\supp \tilde{\varphi}_n = n^2-2n \to \infty$ as $n \to \infty$.   
This, in particular, ensures that 
$\psi_n \in W_0^{1,p}(\Omega_0\setminus K)$ 
for all sufficiently large~$n$.
Let us abbreviate
and 
$
  \|\cdot\|_{\infty} := \|\cdot\|_{L^\infty(\R)}
$
and
$
  \|\cdot\|_{n,\infty} := \|\cdot\|_{L^\infty(\supp \tilde{\varphi}_n)}
$.
Similarly as in the proof of Proposition~\ref{Prop.straight},
we use Lemma~\ref{Lem.crucial} (twice) to obtain
$$
\begin{aligned}
  Q[\psi_n] 
  &\leq  
  \alpha^{p/2} \int_{\Omega_0} 
  \frac{\tilde\alpha^p \, |\partial_s\psi_n|^p 
  + \tilde\beta^p \, |f_\mu\partial_{t_\mu}\psi_n|^p}{f^p} 
  \, f \, ds \, dt
  +  \beta^{p/2} \int_{\Omega_0} 
  |\nabla_{\!t} \psi_n|^p \, f \, ds \, dt 
  \\
  &\leq
  \frac{\alpha^{p/2} }{(1-a\,\|\kappa\|_{n,\infty})^{p-1}}  
  \int_{\Omega_0} 
  \left(
  \tilde\alpha^p \, |\partial_s\psi_n|^p
  + \tilde\beta^p \, |f_\mu\partial_{t_\mu}\psi_n|^p
  \right) \, ds \, dt
  \\
  & \quad
  +  \beta^{p/2} \, (1+a\,\|\kappa\|_{n,\infty})  \int_{\Omega_0} 
  |\nabla_{\!t} \psi_n|^p \,  ds \, dt 
  \\
  &= \frac{\alpha^{p/2} \, \tilde\alpha^p}{(1-a\,\|\kappa\|_{n,\infty})^{p-1}} 
  \int_{\R} |\varphi_n'|^p \, ds
  + \frac{\alpha^{p/2} \, \tilde\beta^p}{(1-a\,\|\kappa\|_{n,\infty})^{p-1}} 
  \int_{\R} |\varphi_n|^p \, r \, ds
  \\
  & \quad
  +  \beta^{p/2} \, (1+a\,\|\kappa\|_{n,\infty}) 
  \, \lambda_1(\omega) \int_{\R} 
  |\varphi_n|^p \, ds  
  \\
  &\leq \frac{\alpha^{p/2} \, \tilde\alpha^p}{(1-a\,\|\kappa\|_{\infty})^{p-1}} 
  \int_{\R} |\varphi_n'|^p \, ds
  + \frac{\alpha^{p/2} \, \tilde\beta^p \, \|r\|_{n,\infty}}
  {(1-a\,\|\kappa\|_{n,\infty})^{p-1}} 
  \int_{\R} |\varphi_n|^p \, ds
  \\
  & \quad
  +  \beta^{p/2} \, (1+a\,\|\kappa\|_{n,\infty}) 
  \, \lambda_1(\omega) \int_{\R} 
  |\varphi_n|^p \, ds  
  \,.
\end{aligned}
$$
At the same time,
$$
  \|\psi_n\|^p \geq  
  (1-a\,\|\kappa\|_{n,\infty})
  \int_{\R} |\varphi_n|^p \, ds  
  \geq (1-a\,\|\kappa\|_{\infty})
  \int_{\R} |\varphi_n|^p \, ds  
  \,.
$$
Consequently,
$$
  \frac{Q[\psi_n]}{\|\psi_n\|^p} \leq 
  \frac{\alpha^{p/2} \, \tilde\alpha^p}{(1-a\,\|\kappa\|_{\infty})^{p}} 
  \, \frac{\displaystyle\int_{\R} |\varphi_n'|^p \, ds}
  {\displaystyle\int_{\R} |\varphi_n|^p \, ds}
  + \frac{\alpha^{p/2} \, \tilde\beta^p \, \|r\|_{n,\infty}}
  {(1-a\,\|\kappa\|_{n,\infty})^{p}} 
  + \beta^{p/2} \, \frac{1+a\,\|\kappa\|_{n,\infty}}{1-a\,\|\kappa\|_{n,\infty}}  
  \, \lambda_1(\omega) 
  \,.
$$
By~\eqref{fnp} and the fact that $\kappa(s) \to 0$ 
and $r(s) \to 0$ as $s \to \infty$, 
it follows that 
$$
  \lim_{n \to \infty} \frac{Q[\psi_n]}{\|\psi_n\|^p}
  \leq \beta^{p/2} \, \lambda_1(\omega) 
  \,,
$$
where $\beta > 1$ can be made arbitrarily close to~$1$.
In summary, given any $K\Subset\Omega_0$ 
and any $\eps>0$,
we have proved that 
$$
  \lambda_1(\Omega_0 \setminus K,g) \leq \lambda_1(\omega) + \eps
  \,.
$$
Consequently, $\lambda_\infty(\Omega_{\kappa,R}) \leq \lambda_1(\omega) + \eps$.
Since~$\eps$ can be made arbitrarily small, 
the desired claim follows.
\end{proof}

As a consequence of Propositions~\ref{Prop.lower} and~\ref{Prop.upper},
we get Theorem~\ref{Thm.ess}.
In particular, the second equality of~\eqref{straight} follows, too.

\begin{rem}
It is clear from the proof of Proposition~\ref{Prop.upper} 
that its simple modification enables one to assume that 
$\kappa(s) \to 0$ and $r(s) \to 0$ as $s \to -\infty$
to get the upper bound  
$\lambda_\infty(\Omega_{\kappa,R}) \leq \lambda_1(\omega)$. 
On the other hand, the lower bound 
$\lambda_\infty(\Omega_{\kappa,R}) \geq \lambda_1(\omega)$
requires that $\kappa(s) \to 0$ both as~$s \to \pm\infty$
(but no condition on~$r$ is needed). 

Inspired by \cite[Ex.~5.1]{kk5},
untwisted periodically non-trivially bent tubes are an example when 
$\lambda_\infty(\Omega_{\kappa,I}) < \lambda_1(\omega)$. 
At the same time, unbent periodically non-trivially twisted tubes satisfy
$\lambda_\infty(\Omega_{0,R}) > \lambda_1(\omega)$.
What is more, if $d=3$ and $p=2$, it is known~\cite{K11}
that $\lambda_\infty(\Omega_{0,R})=\infty$ 
whenever $0 \not\in\omega$ and $r(s) \to \infty$ as $|s| \to \infty$
(an extension of this result to higher dimensions $d \geq 4$
and/or arbitrary $p\in (1,\infty)$ represents
and interesting open problem).
\end{rem}

\section{Bent tubes}\label{bent3}
In this section we establish Theorem~\ref{Thm.disc}
dealing with untwisted bent tubes.
The proof has roots in the original variational idea of~\cite{GJ},
but we rather follow the rigorous implementations 
due to \cite{kk5,cdfk,K-Padova}.
The nonlinear case $p \not= 2$ requires 
technically non-trivial modifications.
What is more, the lack of an integration-by-parts argument
(cf.~Remark~\ref{Rem.circular}) leads us 
to assume that~$\omega$ is circular.

The weak formulation of the eigenvalue equation 
$-\Delta_p\phi_1 = \lambda_1(\omega) |\phi_1|^{p-2} \phi_1$ in~$\omega$,
subject to Dirichlet boundary conditions $\phi_1=0$ on~$\partial\omega$, 
reads
\begin{equation}\label{weak} 
  \forall \varphi \in W_0^{1,p}(\omega)
  \,, \qquad
  \int_{\omega} 
  |\nabla\phi_1|^{p-2} \,
  \nabla\phi_1 \cdot \nabla\varphi
  \, dt
  = \lambda_1(\omega) \int_{\omega} 
  |\phi_1|^{p-2} \,
  \phi_1 \varphi
  \, dt
  \,.
\end{equation}
Unless $p=2$, the eigenvalue $\lambda_1(\omega)$ 
and its associated eigenfunction~$\phi_1$ 
are not known explicitly, even for balls 
(not even in one dimension \cite{Boulton}). 
However, by the positivity and uniqueness of~$\phi_1$,
it is possible to conclude that~$\phi_1$ 
is radially symmetric in balls and spherical shells 
(see, e.g., \cite{Bhattacharya,Nazarov,Ercole,
Bobkov-Drabek_2017,Anoop-Bobkov-Sasi_2018}).
Using this observation, one immediately concludes with 
the following properties that we employ in the proof.
\begin{lem}\label{Lem.symmetry}
If $\omega$ is circular, then
$$
  \int_\omega|\phi_1(t)|^p \, t \, dt = 0
  = \int_\omega |\nabla\phi_1(t)|^p \, t \, dt
  \,.
$$
\end{lem}

Now we are ready to establish Theorem~\ref{Thm.disc}. 
Without loss of generality, we may take $R=I$.

\begin{proof}[Proof of Theorem~\ref{Thm.disc}]
By the definition of $\lambda_1(\Omega_{\kappa,I})$ given in \eqref{ray.bis},
the claim is equivalent to the existence of a (trial) function 
$\psi\in W_0^{1,p}(\Omega_0)$ for which 
\begin{equation}\label{strategy}
  Q_1[\psi] := Q[\psi] - \lambda_1(\omega) \, \|\psi\|^p < 0
  \,.
\end{equation}

The first step consists in 
taking a regularisation of $(s,t) \mapsto \phi_1(t)$,
where~$\phi_1$ is the normalised eigenfunction of 
the Dirichlet $p$-Laplacian in~$\omega$.
More specifically, as in the proof of Proposition~\ref{Prop.critical},
we define 
$$
  \psi_n(s,t) := \varphi_n(s) \phi_1(t)
  \,,
$$
where~$\varphi_n$ is given by~\eqref{one}.  
By virtue of Lemma~\ref{Lem.symmetry}
and the variational definition of $\lambda_1(\omega)$,
one has
\begin{equation}\label{circular}
\begin{aligned}
  \lambda_1(\omega) \, \|\psi_n\|^p 
  &= \lambda_1(\omega) \int_{\Omega_0} 
  |\psi_n|^p \, f \, ds \, dt 
  = \lambda_1(\omega) \int_\R |\varphi_n|^p 
  \int_{\omega} 
  |\phi_1|^p \, f \, dt 
  \, ds
  \\
  &= \lambda_1(\omega) \int_\R |\varphi_n|^p 
  \int_{\omega} 
  |\phi_1|^p \, dt 
  \, ds
  = \int_\R |\varphi_n|^p 
  \int_{\omega} 
  |\nabla_{\!t} \phi_1|^p \, dt 
  \, ds
  \\
  &= \int_\R |\varphi_n|^p 
  \int_{\omega} 
  |\nabla_{\!t} \phi_1|^p \, f \, dt 
  \, ds 
  = \int_{\Omega_0} 
  |\nabla_{\!t}\psi_n|^p \, f \, ds \, dt 
  = \|\nabla_{\!t}\psi_n\|^p
  \,.
\end{aligned}
\end{equation}
Consequently,
$$
  Q_1[\psi_n] = \int_{\Omega_0} 
  \left(
  \left|\frac{\partial_s\psi_n}{f}\right|^2
  + |\nabla_{\!t}\psi_n|^2 
  \right)^{p/2}
  f \, ds \, dt
  -  \int_{\Omega_0} 
  |\nabla_{\!t}\psi_n|^p \, f \, ds \, dt 
  \,.
$$ 
As in the proof of Proposition~\ref{Prop.critical},
we distinguish two cases.
If $p \leq 2$, then
$$
  Q_1[\psi_n] 
  \leq \int_{\Omega_0} 
  \left|\frac{\partial_s\psi_n}{f}\right|^p
  \, f \, ds \, dt 
  \leq \frac{1}{(1 - a\,\|\kappa\|_{L^\infty(\R)})^{p-1}} 
  \int_{\R} |\varphi_n'|^p \, ds 
  \xrightarrow[n \to \infty]{} 0
  \,,
$$ 
where the convergence holds due to~\eqref{fnp}.
If $p>2$, Lemma~\ref{Lem.crucial} yields
$$
\begin{aligned}
  Q_1[\psi_n] 
  &\leq  
  \alpha^{p/2} \int_{\Omega_0} 
  \left|\frac{\partial_s\psi_n}{f}\right|^p \, f \, ds \, dt
  +  (\beta^{p/2}-1)  \int_{\Omega_0} 
  |\nabla_{\!t}\psi_n|^p \, f \, ds \, dt 
  \\
  &\leq \frac{\alpha^{p/2}}{(1 - a\,\|\kappa\|_{L^\infty(\R)})^{p-1}} 
  \int_{\R} |\varphi_n'|^p \, ds 
  +  (\beta^{p/2}-1) \, \lambda_1(\omega) \int_{\R} |\varphi_n|^p \, ds 
  \\
  &\leq \frac{2\alpha^{p/2}}{(1 - a\,\|\kappa\|_{L^\infty(\R)})^{p-1}} 
  \frac{1}{n^{p-1}} 
  +  (\beta^{p/2}-1) \, \lambda_1(\omega) \, 4 n
  \,,
\end{aligned} 
$$
where the last inequality employs~\eqref{fnp} 
and the fact that the function~$\varphi_n$ satisfies
$0 \leq \varphi_n \leq 1$ on $(-2n,2n)$
and that it is zero elsewhere.
Choosing $n$-dependent~$\beta$ as in~\eqref{depend},
it is straightforward to check that 
$Q_1[\psi_n] \to 0$ as $n \to \infty$ 
in this case as well.
In summary, with the choice~$\psi_n$, 
we have achieved an asymptotic equality instead 
of the strict inequality in~\eqref{strategy}.

In the second step, we perturb~$\psi_n$ in such a way that 
the strict inequality is achieved in~\eqref{strategy}.
We define
$$
  \psi_{n,\eps} := \psi_n + \eps \, \phi
  \qquad \mbox{with} \qquad
  \phi(s,t) := j(s) \, \xi(t) \, \phi_1(t)
  \,,
$$
where $\eps \in \R$ and
$j \in C_0^\infty(\R)$ and $\xi \in L^\infty(\omega)$
are arbitrary real-valued functions. 
We always consider~$n$ so large that $\varphi_n=1$ on the support of~$j$. 
Moreover, we always assume $|\eps| \leq \eps_0$,
where $\eps_0 \in \R$ is so small that 
$
  \eps_0 \, \|j\|_{L^\infty(\R)} \, \|\xi\|_{L^\infty(\omega)}
  < 1
$;
this ensures that $\psi_{n,\eps}(s,t) = \phi_1(t) (1 + \eps j(s) \xi(t))$ 
is positive whenever $s \in \supp j$.  
Let us write 
$$
  Q_1[\psi_{n,\eps}] 
  = I_1(\eps) - \lambda_1(\omega) \, I_2(\eps) =: I(\eps)  
  \qquad \mbox{with} \qquad
\begin{aligned}
  I_1(\eps) 
  &:= Q[\psi_{n,\eps}] \,, 
  \\
  I_2(\eps) 
  &:= \|\psi_{n,\eps}\|^p \,.
\end{aligned} 
$$
Our strategy is to employ the Taylor expansion 
\begin{equation}\label{Taylor} 
  I(\eps) = I(0) + I'(0) \, \eps + o(\eps)
  \qquad \mbox{as} \qquad
  \eps \to 0 \,,
\end{equation}
where the remainder~$o(\eps)$ is of Peano type. 

Let us start with the derivative of the simpler integral:
$$
  I_2'(\eps) 
  = p \int_{\Omega_0} |\psi_{n,\eps}|^{p-1} \, \phi \, f \, ds \, dt 
  \,.
$$ 
First of all, observe that it is in fact independent of~$n$,
since $\psi_{n,\eps}(s,t) = \phi_1(t) + \eps \phi(s,t)$ 
whenever $s \in \supp j$. 
Still, we need to argue that the interchange of the derivative 
with respect to~$\eps$ and the integration is justified.
This follows from the $\eps$-independent bound 
$$ 
  \big| |\psi_{n,\eps}(s,t)|^{p-1} \, \phi(s,t) \big|
  \leq  \|j\|_{L^\infty(\R)} \, \|\xi\|_{L^\infty(\omega)} 
  \, (1 + \eps_0 \, \|j\|_{L^\infty(\R)} \, \|\xi\|_{L^\infty(\omega)})^{p-1}
  \, |\phi_1(t)|^p  
  \,,
$$
which is integrable over $(s,t) \in \Omega_0' := \supp j \times \omega$
(the Jacobian~$f$ is irrelevant in view of~\eqref{est}).

As for the derivative of the more complicated integral,
we find 
$$
  I_1'(\eps) = \int_{\Omega_0'} F \, f \, ds \, dt
$$
with
$$
  F
  := p 
  \left(
  \left|\frac{\partial_s\psi_{n,\eps}}{f}\right|^2
  + |\nabla_{\!t}\psi_{n,\eps}|^2 
  \right)^{p/2-1}
  \left(
  \frac{\partial_s\psi_{n,\eps} \, \partial_s\phi}{f^2} 
  + \nabla_{\!t}\psi_{n,\eps} \cdot \nabla_{\!t}\phi
  \right)
  \,.
$$ 
Here the function~$F$ should be interpreted as zero
at the points where both 
$\partial_s\psi_{n,\eps}$ and $\nabla_{\!t}\psi_{n,\eps}$
equal zero.
Again, $I_1'(\eps)$ is actually independent of~$n$.
Indeed, 
$\psi_{n,\eps} = \phi_1 + \eps \phi$ 
and $\partial_s\psi_{n,\eps} = \eps \partial_s\phi$ 
on $\Omega_0'$.
Using the Schwarz inequality, we get the bound
$$ 
  |F|
  \leq p 
  \left(
  \left|\frac{\eps\partial_s\phi}{f}\right|^2
  + |\nabla_{\!t} (\phi_1+\eps\phi)|^2 
  \right)^{(p-1)/2}
  \left(
  \left|\frac{\partial_s\phi}{f}\right|^2
  + |\nabla_{\!t}\phi|^2
  \right)^{1/2}
  \,.
$$
Similarly as above, using additionally that not only $\phi_1$
but also $\nabla_{\!t}\phi_1$ belongs to $L^p(\omega)$,
it is straightforward to estimate this bound
by an $\eps$-independent bound integrable in~$\Omega_0'$.

Let us look at the first variation
$$
  I'(0) = p \int_{\Omega_0} 
  |\nabla_{\!t}\phi_1|^{p-2} \,
  \nabla_{\!t}\phi_1 \cdot \nabla_{\!t}\phi
  \, f \, ds \, dt
  - p \, \lambda_1(\omega) \int_{\Omega_0} 
  |\phi_1|^{p-2} \,
  \phi_1 \phi
  \, f \, ds \, dt
  \,.
$$ 
Choosing $\varphi := \phi(s,\cdot) f(s,\cdot)$ 
for the test function in~\eqref{weak} 
and the Fubini theorem yield
$$
\begin{aligned}
  I'(0) 
  &= p \int_{\Omega_0} 
  |\nabla_{\!t}\phi_1|^{p-2} \,
  \phi \, \nabla_{\!t}\phi_1 \cdot \nabla_{\!t} f
  \, ds \, dt
  \\
  &= -p 
  \int_{\omega} 
  |\nabla_{\!t}\phi_1(t)|^{p-2} \, \xi(t) \, \phi_1(t) \,
  k \cdot \nabla_{\!t} \phi_1(t) \, dt
\end{aligned}  
$$
where $k := (k_1,\dots,k_{d-1})$ is a constant vector 
composed of $k_i := \int_{\R} j(s) \, \kappa_i(s) \, ds$. 
By choosing the support of~$j$ on an interval where $\kappa \not= 0$,  
the vector~$k$ can be chosen to be non-zero. 
We claim that there is a choice of~$\xi$ which guarantees
that $I'(0) \not= 0$. By contradiction, let us assume that 
$I'(0) = 0$ for any choice of~$\xi$.
Then, recalling that~$\phi_1$ is positive,
necessarily $k \cdot \nabla_{\!t} \phi_1 = 0$ in~$\omega$.
Since~$\phi_1$ is radial, it follows that $k \cdot t = 0$   
for every $t \in \omega$, which is a contradiction.

In summary, we have established~\eqref{Taylor} with 
$I(0) = Q_1[\psi_n] \to 0$ as $n \to \infty$, 
$n$-independent remainder $o(\eps)$
and $n$-independent non-zero $I'(0)$.
By choosing~$\eps$ sufficiently small and of suitable sign,
it is possible to ensure that $I'(0) \, \eps + o(\eps) < 0$.
Then we choose~$n$ so large that also $I(\eps) < 0$. 
\end{proof}
\begin{rem}\label{Rem.circular}
Unfortunately, we do not know how to get rid of the hypothesis
that~$\omega$ is circular.
This assumption was employed twice in the proof.
The usage in the argument showing that~$I'(0)$ is non-zero 
can be avoided by noticing that, in general, 
the property that a directional derivative of~$\phi_1$ 
vanishes identically is incompatible with the Dirichlet boundary conditions.
The second circumstance was the usage of the identity
\begin{equation}\label{problem1} 
  \int_\omega |\nabla \phi_1(t)|^p \, t \, dt 
  = \lambda_1(\omega) \int_\omega |\phi_1(t)|^p \, t \, dt 
\end{equation} 
in the first part of the proof (cf.~\eqref{circular}).
This trivially holds for circular domains 
due to the symmetry (cf.~Lemma~\ref{Lem.symmetry}).
By using the test function $\varphi(t):=t \, \phi_1(t)$ in~\eqref{weak},
the identity~\eqref{problem1} is equivalent to the property
\begin{equation}\label{problem2} 
  \int_\omega |\nabla \phi_1(t)|^{p-2} \, \nabla|\phi_1|^2 \, dt 
  = 0
  \,.
\end{equation} 
Integrating by parts, this identity obviously holds 
for any domain~$\omega$ whenever $p=2$. 
For arbitrary $p \in (1,\infty)$,
instead of assuming that~$\omega$ is circular,
the present proof works
(and therefore conclusions of Theorem~\ref{Thm.disc}
and Corollary~\ref{Corol.disc} hold) 
for any domain~$\omega$ satisfying~\eqref{problem2}. 
\end{rem}
\begin{rem} 
Since~$\omega$ is assumed to be circular in Theorem~\ref{Thm.disc},
we could have stated an illusively more general claim 
$\lambda_1(\Omega_{\kappa,R}) < \lambda_1(\omega)$
irrespectively of the choice of~$R$. 
It is just because $\Omega_{\kappa,R}=\Omega_{\kappa,I}$.
\end{rem}

\section{Twisted tubes}\label{Sec.twist}
In this section we establish Theorem~\ref{Thm.Hardy} 
dealing with unbent twisted tubes. 
We therefore assume that $\kappa = 0$ (so $f=1$).
Under the hypothesis~\eqref{twist} that~$\Omega_{0,R}$ 
is non-trivially twisted, 
our goal is to show that there exists  
a positive continuous function $\varrho:\Omega_0 \to \R$
such that the Hardy inequality
\begin{equation}\label{Hardy}
  \forall \psi \in W_0^{1,p}(\Omega_0,g) \,, \qquad
  Q[\psi] - \lambda_1(\omega) \, \|\psi\|^p
  \geq \int_{\Omega_0} \varrho \, |\psi|^p \, ds \, dt
\end{equation}
holds. 
Note that~\eqref{twist} necessarily requires that $d \geq 3$,
because the only orthogonal matrix of dimension~$1$ 
is the scalar identity 
(there is no twisting for two-dimensional strips).  

Following the approach of \cite{K6-with-erratum,KZ1},  
we define 
\begin{equation}\label{ray.compact}
  \lambda_1^N(\Omega_0^l) :=
  \inf_{\stackrel[\psi \not= 0]{}{\psi\in W_0^{1,p}(\Omega_0)}} 
  \frac{\displaystyle Q^l[\psi]}
  {\displaystyle \|\psi\|^{l,p}}
\end{equation}
for any $\Omega_0^l := (-l,l) \times \omega$ with $l>0$, 
where
$$
\begin{aligned}
  Q^l[\psi] &:=
  \int_{\Omega_0^l} 
  \left(
  \big|\big(\partial_s-f_\mu\partial_{t_\mu}\big)\psi\big|^2 
  + |\nabla_{\!t} \psi|^2
  \right)^{p/2} \, ds \, dt 
  \,,
  \\
  \|\psi\|^{l,p} &:= \int_{\Omega_0^l}  |\psi|^p \, ds \, dt \,.
\end{aligned}  
$$
The minimisation over  $\psi\in W_0^{1,p}(\Omega_0)$ in~\eqref{ray.compact}
precisely means that~$\psi$ is a restriction of a function 
from $W_0^{1,p}(\Omega_0)$ to $\Omega_0^l$.
We use the superscript~$N$ to point out that the minimiser~$\psi$ 
of~\eqref{ray.compact} satisfies Neumann boundary conditions 
$(\partial_s-f_\mu\partial_{t_\mu}\big)\psi=0$
on $\{\pm l\} \times \omega$,
but we shall not use this fact.

By~\eqref{Poincare}, 
$\lambda_1^N(\Omega_0^l) \geq \lambda_1(\omega)$.
At the same time, by choosing the trial function $\psi(s,t) := \phi_1(t)$
in~\eqref{ray.compact},
it follows that $\lambda_1^N(\Omega_0^l) = \lambda_1(\omega)$ if 
$f_\mu\partial_{t_\mu} \phi_1 = 0$ in~$\Omega_0^l$.
The converse result is non-trivial.
\begin{lem}\label{Lem.lambda}
One has
$$
  \lambda_1^N(\Omega_0^l) > \lambda_1(\omega)
  \qquad \Longleftrightarrow \qquad
  f_\mu\partial_{t_\mu} \phi_1 \not= 0 
  \quad \mbox{in} \quad \Omega_0^l
  \,.
$$
\end{lem}
\begin{proof}
To prove the remaining implication, let us assume that 
$f_\mu\partial_{t_\mu} \phi_1 \not= 0$ 
but $\lambda_1^N(\Omega_0^l) = \lambda_1(\omega)$.
By compactness, the infimum~\eqref{ray.compact}
is indeed achieved by a function 
$\psi_1 \in W_0^{1,p}(\Omega_0) \upharpoonright \Omega_0^l$.
Then 
$$
  \int_{\Omega_0^l} 
  \big|\big(\partial_s-f_\mu\partial_{t_\mu}\big)\psi_1\big|^p
   \, ds \, dt = 0
  \qquad \mbox{and} \qquad
  \int_{\Omega_0^l} 
  |\nabla_{\!t} \psi_1|^p
  \, ds \, dt 
  - \lambda_1(\omega) \int_{\Omega_0^l} |\psi_1|^p = 0
  \,.
$$  
By~\cite{Kawohl-Lindqvist},  
$\lambda_1^N(\Omega_0^l)$ is simple.
Consequently, the second identity implies  
that there exists a function $\varphi \in W^{1,p}((-l,l))$
such that $\psi_1(s,t) = \varphi(s) \phi_1(t)$.  
Substituting this result 
into the first identity, we get 
\begin{equation}\label{divergence}
  0 = \varphi' \phi- \varphi  f_\mu \partial_{t_\mu}\phi_1 
  = \varphi' \phi - \varphi  \partial_{t_\mu} (f_\mu \phi_1)
  = \divergence (\varphi \phi_1, - \varphi f_\mu \phi_1)
  \,,
\end{equation}
where the second equality follows by the orthogonality of~$R$.
By the divergence theorem, it is possible to conclude 
that $\varphi$~is constant. 
Substituting this result back to~\eqref{divergence},
it follows that
$f_\mu\partial_{t_\mu} \phi_1 = 0$ in $\Omega_0^l$, a contradiction.
\end{proof}

Now we are ready to establish Theorem~\ref{Thm.Hardy}.
\begin{proof}[Proof of Theorem~\ref{Thm.Hardy}]
Definition~\eqref{ray.compact} implies 
$$
 \forall \psi \in W_0^{1,p}(\Omega_0,g) \,, \qquad
  Q[\psi] - \lambda_1(\omega) \, \|\psi\|^p
  \geq 
  c_l
  \int_{\Omega_0^l} |\psi|^p 
$$
for every positive~$l$, 
where $c_l := \lambda_1^N(\Omega_0^l) - \lambda_1(\omega)$.
By Lemma~\ref{Lem.lambda} and hypothesis~\eqref{twist},
the constant~$c_l$ is positive
for all sufficiently large~$l$. 
This establishes a ``local'' Hardy inequality with
$$
  \varrho_l := c_l \, \chi_{\Omega_0^l}
  \,.
$$ 
We call it local, because the weight~$\varrho_l$ is not positive,
albeit it is non-negative and non-trivial.
However, there exists a general procedure how to deduce
a ``global'' Hardy inequality (i.e., with a positive~$\varrho$)
from the local one. 
It is based on a standard argument of partition of unity 
subordinated to a finitely local covering, 
see \cite[Lem.~3.1]{Pinchover-Tintarev_2007}.
In detail, given any natural number $j \geq 1$, let us write
$$
\begin{aligned}
  2^{-j} \big( Q[\psi] - \lambda_1(\omega) \, \|\psi\|^p \big) 
  &\geq 2^{-j} \, c_{l+j} \int_{\Omega_0^{l+j}} |\psi|^p \, ds \, dt
  \\
  &\geq 2^{-j} \, \min\{c_{l+j},1\} \int_{\Omega_0} \chi_{\Omega_0^{l+j}} 
  \, |\psi|^p \, ds \, dt
  \,.
\end{aligned}  
$$
Summing over all $j \geq 1$ and interchanging 
the order of summation and integration, 
we get~\eqref{Hardy} with 
$$
  \varrho(s,t) := \sum_{j=1}^\infty 
  2^{-j} \, \min\{c_{l+j},1\} \, \chi_{[-(l+j),l+j]}(s)
  \,.
$$   
Since this Hardy weight is independent of~$t$,
one gets~\eqref{Hardy.intro} with $\rho(s,t) := \varrho(s,t)$. 
\end{proof}

Finally, let us comment on hypothesis~\eqref{twist}.
\begin{rem}[$d=3$]
In three dimensions, one has a convenient parameterisation
$$
  R(s) =
  \begin{pmatrix}
    \cos\theta(s) & -\sin\theta(s) \\
    \sin\theta(s) & \cos\theta(s)
  \end{pmatrix}
  ,
$$
where $\theta:\R\to\R$ is a differentiable function
with locally bounded derivative. 
Then condition~\eqref{twist} 
is equivalent to a simultaneous validity 
of the following two requirements:
$$
  \theta' \not = 0
  \qquad \mbox{and} \qquad 
  \omega \mbox{ is not circular}
  \,.
$$
This is clear from the identity
$
  f_\mu \partial_\mu = \theta' \, \partial_\tau
$,
where $\partial_\tau := t_2 \partial_{t_1} -t_1 \partial_{t_2}$
is the angular derivative.
\end{rem}
\begin{rem}[$d \geq 4$]
In higher dimensions, the situation is more complicated
because we cannot separate the ``longitudinal'' and ``transverse''
variables from the condition (which is natural in view of
a more complicated structure of rotations in the higher dimensions).
Anyway, we have the following sufficient condition:
$$
\left.
\begin{aligned}
  &R'\not=0
  \\
  &\forall \mbox{tangential } \sigma\in\R^{d-1}, \
  \sigma\not=0,
  \qquad
  \partial_\sigma \phi_1 \not=0
\end{aligned}
\right\}
  \qquad \Longrightarrow \qquad
  \eqref{twist} \mbox{ holds}
  \,.
$$
(By a ``tangential'' vector in~$\R^{d-1}$ we mean
any vector perpendicular to the radial vector $t \in \R^{d-1}$.)
The implication is clear from the identity
(employing the orthogonality condition $RR^T = I$)
$$
  f_\mu t_\mu = t_\alpha R_{\alpha\beta}' R_{\mu\beta} t_\mu
  = - t_\alpha R_{\alpha\beta} R_{\mu\beta}' t_\mu
  = - t_\mu f_\mu
  \,,
$$
which shows that $(f_1(s),\dots,f_{d-1}(s))$ is tangential
for every $s \in \R$.
In particular, another sufficient condition follows:
$$
\left.
\begin{aligned}
  & R'\not=0
  \\
  & 0 \not\in \omega
\end{aligned}
\right\}
  \qquad \Longrightarrow \qquad
  \eqref{twist} \mbox{ holds}
  \,.
$$
\end{rem}

\section*{Acknowledgments}
D.K.\ was supported
by the EXPRO grant No.~20-17749X
of the Czech Science Foundation.
This work has been partially carried out during a stay of L.B. in Prague at the \emph{Department of Mathematics}, 
\emph{Faculty of Nuclear Sciences and Physical Engineering}, 
\emph{Czech Technical University in Prague}. 
She would like to express her deep gratitude to this prestigious institution for its support and warm hospitality.
L.B.\ is member of the {\em Gruppo Nazionale per l'Analisi Ma\-te\-ma\-ti\-ca, la Probabilit\`a e le loro Applicazioni}
(GNAMPA) of the {\em Istituto Nazionale di Alta Matematica} (INdAM).
L.B. was partially supported by National Science Centre, Poland (Grant No. 2020/37/B/ST1/02742), by INdAM-GNAMPA Project 2023 titled {\em Problemi ellittici e parabolici con termini di reazione singolari e convettivi} (E53C22001930001) and by the IMAG-Maria de Maeztu Excellence Grant CEX2020-001105-M funded by MICINN/AEI.

{\bf Data Availability:} The manuscript has no associated data.

{\bf Conflict of interest:} On behalf of all authors, the corresponding author states that there is no conflict of interest.

%
\bibliography{pbib}
\bibliographystyle{amsplain}

\end{document}